\documentclass[12pt,a4paper,reqno]{amsart}
\usepackage{amsmath,amsthm}
\usepackage{amsfonts}
\usepackage{amssymb}
\usepackage{tikz-cd}
\allowdisplaybreaks

\newcommand{\be}{\begin{equation}}
\newcommand{\ef}{\end{equation}}
\chardef\bslash=`\\ 

\newtheorem*{thm*}{Theorem}

\theoremstyle{definition}

\newtheorem*{remark*}{Remarks}
\newtheorem*{defn*}{Definition}

\theoremstyle{remark}

\newcommand{\wt}{\widetilde}
\newcommand{\wh}{\widehat}
\newcommand{\fc}{\frac}

\newcommand{\iy}{\infty}
 \renewcommand{\sectionmark}[1]{}
\renewcommand{\Im}{\operatorname{Im}}

\newcommand{\const}{\operatorname{const}}

 \usepackage{amsfonts}

\newcommand{\dl}{\delta}
\newcommand{\D}{\mathbb{D}}
\newcommand{\om}{\omega}
\newcommand{\z}{\zeta}
\newcommand{\ov}{\overline}
\newcommand{\vp}{\varphi}
\newcommand{\hC}{\wh{\mathbb{C}}}
\newcommand{\C}{\mathbb{C}}

\newcommand{\B}{\mathbf{B}}
\newcommand{\T}{\mathbf{T}}
\newcommand{\Belt}{\operatorname{Belt}}

\newcommand{\dist}{\operatorname{dist}}

\newcommand{\Om} {\Omega}
\newcommand{\vk} {\varkappa}
\newcommand{\x} {\mathbf x}
\renewcommand{\a} {\alpha}

\begin{document}

\title{All Teichm\"{u}ller spaces are not starlike}
\author{Samuel L. Krushkal}

\begin{abstract}
This paper is the final step in solving the problem of starlikeness of Teichm\"{u}ller spaces in Bers' embedding. This step concerns the case of finite dimensional Teichm\"{u}ller spaces $\T(g, n)$ of positive dimension (corresponding to punctured Riemann surfaces of finite conformal type $(g, n)$ with $2g - 2 + n > 0$).
\end{abstract}

\date{\today\hskip4mm ({\tt TnotStar(3).tex})}

\maketitle

\bigskip

{\small {\textbf {2000 Mathematics Subject Classification:} Primary:
30C62, 30F60, 32G15}}

\medskip

\textbf{Key words and phrases:} Teichm\"{u}ller space, holomorphic embedding, Schwarzian derivative, starlike

\bigskip

\markboth{Samuel L. Krushkal}{Teichm\"{u}ller spaces not starlike}
\pagestyle{headings}

\bigskip\noindent
{\bf 1. Preamble and main result}.
The question on shape of holomorphic embeddings of Teichm\"{u}ller  spaces was raised in the book \cite{BK}
and states:

\bigskip\noindent
{\it For an arbitrary finitely or infinitely generated Fuchsian group $\Gamma$ is the Bers embedding of its Teichm\"{u}ller space $\mathbf T(\Gamma)$ starlike?}

\bigskip
In this embedding, the space $\mathbf T(\Gamma)$ is represented as a bounded domain formed by the Schwarzian derivatives
$$
S_w = \Bigl(\frac{w^{\prime\prime}}{w^\prime}\Bigr)^\prime - \frac{1}{2} \Bigl(\frac{w^{\prime\prime}}{w^\prime}\Bigr)^2
$$
of holomorphic univalent functions $w(z)$ in the lower half-plane $U^* = \{z: \Im z < 0\}$)
(or in the disk) admitting quasiconformal extensions to the Riemann sphere $\wh{\C} = \C \cup \{\iy\}$ compatible with the group $\Gamma$ acting on $U^*$.

The first result here was established by the author in 1989 answering negatively for universal Teichm\"{u}ller space.
It was shown in \cite{Kr1} that universal Teichm\"{u}ller space $\mathbf T = \mathbf T(\mathbf 1)$ has
points which cannot be joined to a distinguished point even by curves of a considerably general form, in particular, by polygonal lines with the same finite number of rectilinear segments. The proof relies on the existence of conformally rigid domains established by Thurston in \cite{Th} (see also \cite{As}).

This implies that the space $\T$ is not starlike with respect to any of its points, and there exist points $\varphi \in \mathbf T$ for which the line interval
$\{t \varphi: 0 < t < 1\}$ contains the points from $\mathbf B \setminus \mathbf S$, where
$\mathbf B = \mathbf B(U^*)$ is the Banach space of hyperbolically bounded holomorphic functions in the
half-plane $U^*$ with norm
$\|\vp\|_{\mathbf B} = 4 \sup_{U^*} y^2 |\vp(z)|$
and $\mathbf S$ denotes the set of all Schwarzian derivatives of univalent functions on $U^*$. These points correspond to holomorphic functions on $U^*$ which are only locally univalent.

On this way, it was established in \cite{Kr3} and \cite{To} that also all finite dimensional Teichm\"{u}ller spaces $\mathbf T(\Gamma)$ of high enough dimensions $n \ge n_0$ and the spaces corresponding to Fuchsian groups of second kind. are not starlike.

Recently, it was established in \cite{Kr7} by a different approach that all spaces $\T(g, 0)$ of closed Riemann surfaces of genus $g \ge 2$ also are not starlike. The proof of this is constructive and provides
the surfaces violating starlikeness.

\bigskip
The present paper is the final step in solving this problem. It concerns the remaining case of finite dimensional Teichm\"{u}ller spaces $\T(g, n)$ of positive dimension corresponding to punctured Riemann surfaces of finite conformal type $(g, n)$(surfaces of genus $g$ with $n$ punctures, up to conformal equivalence).
We show that the answer also is negative also in this case:

\bigskip\noindent
{\bf Theorem 1}. {\it All Teichm\"{u}ller spaces $\T(g, n)$ with $2g - 2 + n > 0$
are not starlike in Bers embedding}.

\bigskip
This theorem  includes all spaces $\T(g, n)$ of positive dimension. So, all Teichm\"{u}ller spaces $\T(\Gamma)$ representing hyperbolic Riemann surfaces  are not starlike.

\bigskip
Recall that the antipodal quantity $\chi = - (2g - 2 + n)$ represents the Euler characteristic of surface $X = \D/\Gamma$ (and of its quasiconformal deformations forming the space $\T(g, n)$), and
$$
\dim \T(g, n) = 3g - 3 + n.
$$

\bigskip\noindent
{\bf 2. Two basic lemmas}.
We shall represent the points of $\T(\Gamma)$ by Fuchsian groups $\Gamma$ acting discontinuously on
the disks $\D = \{|z| < 1\}$ and $\D^* = \{z \in \hC: |z| > 1\}$ and use the corresponding space $\B = \B(\D)$ of hyperbolically bounded functions$\vp$ on $\D$ with norm
$$
\|\vp\|_{\B} = \sup_\D (1 - |z|^2)^2 |\vp(z)|.
$$
The similar space on the disk $\D^*$ is formed by holomophic functions in the disk $\D^*$
satisfying $\vp(z) = O(z^{-4})$ near $z = \iy$.

For any Fuchsian group $\Gamma$ of acting on the disk $\D$,
the universal covering $\D \to \D/\Gamma$ naturally generates a canonical embedding of the corresponding Teichm\"{u}ller space $\T(\Gamma)$ into the universal space $\T$. All our groups $\Gamma$ are of the first kind, that means the unit circle $\mathbb S^1 = \partial \D$ is the limit set of $\Gamma$.

Accordingly, we have on each space $\T(\Gamma)$ two Teichm\"{u}ller distances generated by quasiconformal maps: the restricted metric $\wt \tau_{\T(\Gamma)} = \tau_\T|\T(\Gamma)$ and the intrinsic metric $\tau_{\T(\Gamma)}$ generated by maps compatible with the group $\Gamma$.
These metrics are equivalent (see\cite{Le}, Ch. 5; \cite{St}).
Below we shall use on the space $\T(\Gamma)$ the distance $\wt \tau_{\T(\Gamma)}$ coming from the ambient space $\T$  and write, if needed, $(\T(\Gamma), \wt \tau)$.

Recall also that every space $\T(\Gamma)$ of dimension greater than $1$, does not be equivalent holomorphically to a Banach ball, and
$\T(\Gamma) = \T \cap \B(\Gamma)$,
where $\B(\Gamma)$ is the $(3n - 3)$-dimensional linear subspace of $\B$ formed by the $\Gamma$-automorphic forms on $\D$ of weight $- 4$ (quadratic $\Gamma$-differentials), and $\T(\Gamma)$ contains the ball
$\{\|\vp\|_{\B(\Gamma)} < 2\}$, there are the points $\vp \in \T(\Gamma)$ with $\|\vp\|_{\B(\Gamma)} > 2$; cf. \cite{Ea}, \cite{Le}.

\bigskip\noindent
{\bf Lemma 1}. \cite{Kr5} {\it For any rational function $r_n$ with poles of order two on the boundary circle $\mathbb S^1$ of the form
$$
r_n(z) = \sum\limits_1^n \fc{c_j}{(z - a_j)^2} + \sum\limits_1^n \fc{c_j^\prime}{z - a_j}
$$
satisfying $\sum\limits_1^n |c_j| > 0$, we have the equality
 \be\label{1}
\|r_n\|_\B = \limsup\limits_{|z| \to 1} (1 - |z|^2)^2 \ |r_n(z)|.
\end{equation}
}

It follows from (1) that there is a boundary point $z_0 \in \mathbb S^1$ at which the maximal value of $(1 - |z|^2)^2 |r_n(z)|$ on the closed disk $\ov \D$ is attained.

\bigskip
The following lemma involves the Grunsky coefficients of univalent functions in the disk and the extremal Beltrami coefficients of quasiconformal extensions of such functions.

Recall that the Grunsky coefficients of a univalent function $f(z)$ on the disk $\D^*$ are defined from the expansion
 \be\label{2}
\log \fc{f(z) - f(\z)}{z - \z} = - \sum\limits_{m,n=1}^\iy \a_{mn} z^{-m} \z^{-n}, \quad (z, \z) \in (\D^*)^2,
\end{equation}
and by the Grunsky theorem a holomorphic function $f(z) = z + \const + O(z^{-1})$
in a neighborhood of $z = \iy$ can be extended to a univalent holomorphic function on $\D^*$ if and only
the coefficients $\a_{m n}$ satisfy the inequality
$$
\Big\vert \sum\limits_{m,n=1}^\iy \ \sqrt{m n} \ \a_{m n} x_m x_n \Big\vert \le 1,
$$
for any sequence $\mathbf x = (x_n)$ from the unit sphere $S(l^2)$ of the Hilbert space $l^2$ with norm $\|\x\| = (\sum\limits_1^\iy |x_n|^2)^{1/2}$ (here the principal branch of the logarithmic function is chosen); see \cite{Gr}, \cite{Po}. The quantity
$$
\vk(f) = \sup \Big\{\Big\vert \sum\limits_{m,n=1}^\iy \ \sqrt{mn} \ \a_{m n} x_m x_n \Big\vert: \
\mathbf x = (x_n) \in S(l^2) \Big\} \le 1
$$
is called the {\bf Grunsky norm} of $f$. This norm regarded as a function of $S_f$ is a plurisubharmonic function on the space $\T$.

It is majorated by the {\bf Teichm\"{u}ller norm} $k(f^\mu)$ equal to the minimal dilatation   $\inf \|\mu\|_\iy$ of quasiconformal extensions of $f$ to $\D^*$, and the equality $\vk(f)= k(f)$ is valid if and only if the extremal Beltrami coefficient $\mu_0$  in the equivalence class of $f$ (the collection of maps $f^\mu$ extending $f$ to $\hC$) satisfies \cite{Kr2}
$$
\|\mu_0\|_\iy = \sup_{\psi \in A_1^2(\D),\|\psi\|_{A_1} =1}
\Big\vert \iint\limits_\D \mu_0(z) \psi(z) dx dy \Big\vert \quad (z = x + iy \in \D),
$$
where $A_1(\D^*)$ is the space of integrable holomorphic quadratic differentials $\psi(z) dz^2$ on
$\D^*$ and
$$
A_1^2(\D) = \{\psi = \om^2 \in A_1(\D): \ \om \ \ \text{holomorphic \ in} \ \ \D\}
$$
On the other hand, a Beltrami coefficient $\mu_0 \in \Belt(D)_1$ is extremal if and only if
 \be\label{3}
\|\mu_0\|_\iy = \sup_{\|\psi\|_{A_1} =1}
\Big\vert \iint\limits_\D \mu_0(z) \psi(z) dx dy \Big\vert .
\end{equation}

Now, let $L$ be a piecewise $C^{1+}$-smooth positively oriented quasicircle separating the points $0$ and $\iy$, with the interior and exterior domains $D$ and $D^*$ (so $D \ni 0$ and $D^* \ni \iy$),
such the exterior conformal map $f: \ \D^* \to D^*$ has the expansion
$$
f(z) = z + b_0 + b_1 z^{-1} + \dots
$$
(i.e., $f(\iy) =\iy, \ f^\prime(\iy)= 1$, and its Grunsky coefficients are given by (4)).
Assume that $L$ is the union of a circular polygonal line $L^\prime$ with endpoints $\z_1, \ \z_2$ and of complementary $C^{1+}$-smooth arc $L^{\prime\prime}$ so that $L$ is smooth at the points $\z_1, \ \z_2$.

\bigskip\noindent
{\bf Lemma 2}. {\it For any domain $D^*$ satisfying the above assumptions on its boundary and for any $t$
satisfying $|t| \|S_f\|_\B < 2$, the (univalent) solutions $w_t(z) = w(z, t)$ of the Schwarz equation
$$
\Bigl(\frac{w^{\prime\prime}}{w^\prime}\Bigr)^\prime - \frac{1}{2} \Bigl(\frac{w^{\prime\prime}}{w^\prime}\Bigr)^2 = t S_f(z) \quad (z \in \D^*)
$$
with $w_t(\iy) = \iy,\ w_t^\prime(\iy) = 1$ have equal Grunsky and Teichm\"{u}ller norms given by}
 \be\label{4}
\vk(w_t) = k(w_t) = |t| \|S_f\|_\B.
\end{equation}

\noindent
{\bf Proof}. By the Ahlfors-Weill theorem \cite{AW}, the functions $w_t(z)$ are univalent on $\D^*$ and have quasiconformal extensions to $\D$ with harmonic Beltrami coefficients
 \be\label{5}
\nu_{w_t}(z) = - \fc{1}{2} (1 - |z|^2)^2 \ S_{w_t}(1/\ov z) = - \fc{t}{2} (1 - |z|^2)^2 \ S_f(1/\ov z),
\quad |z| <  1.
\end{equation}

Lemma 1 (actually, the equality (1)) allows one to establish by applying to the maps $w_t$ the arguments from  \cite{Kr6}, \cite{Kr7} that the harmonic Beltrami coefficients (5) are extremal in their equivalence classes and satisfy the equalities (4).

\bigskip\noindent
{\bf Remarks}.
Lemma 2 admits a straightforward extension to a more general situation, when the boundary curve $L$ is
the union of the finite collections of arcs of type $L^\prime$ and $L^{\prime\prime}$.
Such a case appears in the proof of Theorem 1.

\bigskip
Note also that the common value of terms in (4) is equal to $1 - \a$, where $\pi \a$ is the least inner angle of the interior domain $D$ and that the maps $w_t$ with $t$ satisfying $|t| \|S_f\|_\B < 2$ do not have the extremal quasiconformal extension across the unit circle $\mathbb S^1$ of Teichm\"{u}ller type

\bigskip\noindent
{\bf 3. Proof of Theorem 1}.
Using the above lemmas, one can obtain the proof of Theorem 1 in the same line as the corresponding result for the spaces of closed Riemann surfaces in \cite{Kr8}.

Consider the space $T(\Gamma) = \T(g, n)$ with $n \ge 1$, formed by punctured Riemann surfaces $X$ of type $(g, n)$ and assume first that  $\T(g, n) \ne \T(1, 1)$, hence $\dim \T(g, n) > 1$.

Fix a base point $X_0$ of the space $\T(g, n) = \T(X_0)$ represented as the quotient space $X_0 = \D^*/\Gamma_0$ with a Fuchsian group $\Gamma_0$ acting discontinuously  on $\D^*$ and $\D$, and take a surface $X \in \T(X_0)$ corresponding to a point $\vp \in \B(\Gamma_0)$ with
 \be\label{6}
\|\vp\|_{\B(\Gamma_0)} > 2.
\end{equation}
This point is represented by a quasifuchsian group
$\Gamma_\mu = (f^\mu)^{-1} \Gamma_0 f^\mu$, whose conjugating quasiconformal automorphism $f^\mu(z)$ of $\hC$ (compatible with $\Gamma_0$) is conformal on $\D^*$.

Consider the Ford fundamental polygon $P(\Gamma_\mu)$ of this group in domain $(D^\mu)^* = f^\mu(\D^*)$ centered at $z = \iy$. It has a finite number of sides (the arcs of the isometric circles of Moebius transformations $\gamma$
generating the group $\Gamma^\mu$); these sides are pairwise $\Gamma_\mu$-equivalent.

In contrast to the case of closed Riemann surfaces, this polygon has a finite number of parabolic vertices $p_l$
corresponding to the punctures of $X_0$. These vertices are located on the boundary curve $f^\mu(\mathbb S^1)$, and the inner angles of $P(\Gamma_\mu)$ at these vertices are equal zero. Thus the boundary of polygon $P(\Gamma_\mu)$
s not a quasiconformal curve, and we must use the appropriately truncated polygons.

Delete from $P(\Gamma_\mu)$ the circular triangles $\Delta(p_l,\dl_l)$ chosen so that one of the vertices of  $\Delta(p_l,\dl_l)$ is the puncture $p_l$, his lateral sides emanate from $p_l$ and are orthogonal to
the third side $\sigma_l$ of this triangle; here
$\delta =\dist (p_l, \sigma_l)$. Set
$$
P_\delta(\Gamma_0) = P(\Gamma_0) \setminus  \bigcup \Delta(p_l,\dl_l),
$$
and replace (smooth) the polygonal lines formed by the sides of $\Delta(p_l,\dl_l)$ opposite to vertices $p_l$ by
a smooth line up to its endpoints.

Now, denoting the elements of $\Gamma_\mu$ by $\gamma_1 = \mathbf 1, \ \gamma_2, \gamma_3, \ \dots $,
numerated so that two successive transformations $\gamma_j, \gamma_{j+1}$ determine the adjacent polygons
$ \gamma_j P(\Gamma_\mu), \  \gamma_{j+1} P(\Gamma_\mu)$ with a common side and choosing a sequence $\delta_n \to 0$, we construct the following increasing collection of curvelinear polygons
 \be\label{7}
\begin{aligned}
\Om_{l,1} &= P_{\dl_1}(\Gamma_\mu), \quad \Om_{l,2} = \gamma_1 P{\dl_2}(\Gamma_\mu) \bigcup \gamma_2 P{\dl_2}(\Gamma_\mu),    \\
\quad \Om_{l,3} &= \gamma_1 P{\dl_3}(\Gamma_\mu) \bigcup \gamma_2 P{\dl_3}(\Gamma_\mu) \bigcup\gamma_3 P{\dl_3}(\Gamma_\mu),  \\
\quad \Om_{l,4} &= \gamma_1 P{\dl_4}(\Gamma_\mu) \bigcup \gamma_2 P{\dl_4}(\Gamma_\mu) \bigcup\gamma_3 P{\dl_4}(\Gamma_\mu) \bigcup \gamma_4 P{\dl_4}(\Gamma_\mu), \dots \ ,
\end{aligned}
\end{equation}
with vertices coming from the initial fundamental polygon $P(\Gamma_\mu)$.
As $l \to 0, \ j \to \iy$, the domains $\Om_{l,j}$  exhaust increasingly the domain $D_\mu^*$.

Take the conformal maps $F_{j,l}$ of the disk $\D^*$ onto these domains and define  their harmonic Beltrami coefficients
 \be\label{8}
t \nu_{F_{j,l}}(z) = - \fc{t}{2} (1 - |z|^2)^2 \ S_{F_{j,l}}(1/\ov z), \quad |z| <  1,
\end{equation}
for $t$ satisfying $|t| \|S_{F_{j,l}}\|_\B < 2$. Similar to Lemma 2, one obtains that all these harmonic coefficients are extremal in their equivalence classes; in addition, the corresponding maps $F_{j,l;t}$ have equal Teichm\"{u}ller and Grunsky norms.

Since the domains (7) exhaust increasingly the Jordan (quasiconformal) domain $D_\mu$, one obtains from the
general Carath\'{e}odory theorem on convergence of conformal maps from the disk onto Jordan domains that the limit function
 \be\label{9}
F(z) = \lim\limits_{l \to 0, j \to \iy} F_{j,l}(z)
\end{equation}
maps conformally the disk $\D^*$ onto the exterior component $D_\mu^*$ of domain of discontinuity of the group $\Gamma_\mu$.
Note that the convergence in (10) is uniform on the disk $\D^*$ in the spherical metric on $\hC$.

Our goal now is to establish that one of the extremal qusiconformal extensions $\wh F^\mu$ of this limit function  inherits the properties of the functions $F_{j,l}$; namely, its harmonic coefficient also is extremal.
This is one of the crucial steps in the proof of Theorem 1.

It follows from (9) that the Beltrami coefficients (8) are convergent on the disk $\D$ to the harmonic
Beltrami coefficient $\nu_{S_F}$ of function $F$; in addition, the general properties of quasiconformal maps imply
 \be\label{10}
\|S_F\|_\B \le \lim\limits_{l \to 0, j \to \iy} \|S_{F_{j,l}}\|_\B.
\end{equation}

We have to establish that in our situation generated by Riemann surfaces of finite type $(g, n)$ we have in (10) the case of equality. This implies that for admissible $|t| > 0$ the harmonic
Beltrami coefficients $\nu_{t S_F}$ also are extremal in their equivalence classes (in the second distance $\wt \tau_{\T(\Gamma)}$ on $\T(\Gamma)$ coming from $\T$).

First observe that normalising the quasiconformal extensions $\wh F(z)$ of univalent functions $F(z) = z + b_0 + b_1 z^{-1}+ \dots$ on $\D^*$ to $\hC$ additionally by $\wh F(0) = 0$, one obtains that their
Teichm\"{u}ller norm $k(F)$ is plurisubharmonic on the universal Teichm\"{u}ller space $\T$ (hence,
also on $\T(\Gamma_\mu)$ in both distances $\tau_{\T(\Gamma)}$ and $\wt \tau_{\T(\Gamma)}$; see, e.g., \cite{Kr4}).

Now we pass to the homotopy functions
$$
F_{j,l,s}(z) = s F_{j,l}(z/s) = z + b_{0,j,l} s + b_{1,j,l} s^2 z^{-1} + \dots
$$
with $|s| \le 1$. The dilatations $k(F_{j,l,s})$ of these functions are circularly symmetric with respect to $s \in \D$ and hence continuous in $|s|$ on $[0, 1]$.

But for any $s$ with $|s| \in (0, 1)$, we have
$$
\|S_{F_{j,s}} - S_{F_s}\|_\B \to 0 \quad \text{as} \ \ j \to \iy,
$$
and the plurisubharmonicity of $\vk(f^\mu)$ on $\T$ implies
$$
\lim\limits_{s \to 1} \vk(F_{j,l,s}) = \vk(F_{j,l})
$$
(cf. \cite{Kr5}, \cite{Ku}), while the properties of extremal quasiconformal maps yield for fixed $l$ and $j$,
$$
\lim\limits_{s \to 1} k(F_{j,l,s}) = k(F_{j,l}).
$$
In view of the relation between the norms $\vk(f)$ and $k(f)$, all this results for $\|S_F\|_\B < 2$ in
the relations
$$
\lim\limits_{s \to 1} k(F_{j,l,s}) = k(F_{j,l}) = \|\nu_{S_{F_{j,l}}}\|_\iy,
$$
which implies the desired equality in (10).

It follows also that all harmonic Beltrami coefficients $\nu_{t S_F}$ with norm less than $1$  must be extremal in their equivalence classes, and therefore, as $\|\nu_{t S_F}\|_\B \to 1$, the corresponding value  $t_0$ must define a boundary point of both spaces $\T(\Gamma_\mu)$ and $\T$ (i.e., $t S_F \to t_0 S_F \in \partial \ T(\Gamma_\mu)$ as $t \to t_0$).

This yields, since $\T(\Gamma)$ and its complementary domain in $\B(\Gamma)$  have a common boundary,  that the subinterval of the ray ${r S_F: \ r > 0}$ between the points  $r_0 S_F $ and $r_{*} S_F$ (the second endpoint corresponds to the surface $D_\mu^*/\Gamma_\mu$ obeying (6)) cannot lie entirely in $\T(\Gamma)$.
This proves the assertion of Theorem 2 for spaces $\T(g, n)$ of dimension greater than $1$.

Now consider the one-dimensional case omitted above. In this case, there exist two
one-dimensional Teichm\"{u}ller spaces $\T(1, 1)$ of punctured tori and $\T(0, 4)$ with quadruples of punctures.
These spaces are canonically isomorphic and either of those is conformally equivalent to the unit disk.

We need for these spaces the following result.

\bigskip\noindent
{\bf Lemma 3}. {\it The space $\T(\Gamma) = \T(1, 1)$ (and similarly $\T(0, 4)$) is broader than the disk } $\{\vp \in \B(\Gamma): \ \|\vp\|_{\B(\Gamma)} < 2\}$.

\bigskip
The fact that the space $\T(\Gamma)$ is not a ball in the space $\B(\Gamma)$ is well-known for spaces of dimension greater than $1$ (being obtained, for example, in \cite{Ea});
\footnote{Note that for $n > 1$, this follows also from the established recently fact that the invariant
Carath\'{e}odory and Kobayashi metrics on $\T(g,n)$ are different.}
it can be derived for the space $\T(0, 4)$ from numerical calculations in \cite{Po}.

For completeness, we provide here another proof for $\T(0, 4)$, which is direct and simpler.
Fix a quadruple $(0, 1, \a, \iy)$ and consider the punctured sphere
$X_0 = \hC \setminus {0, 1, \a,\iy} = \D/\Gamma$ as the base point of the space $\T(0, 4)$.
The Beltrami differentials $\mu(z) d \ov z/dz$ on $X_0$ are lifted to $(-1,1)$-measurable forms $\wh \mu$
on $\D$, compatible with the group $\Gamma$, and being extended by zero to the complementary disk
$\D^*$ determine quasiconformal automorphisms $w^{\wt \mu}$ of $\hC$
conformal on $\D^*$.

Recall that $\T(\Gamma) = \T \bigcap \B(\Gamma)$, where $\B(\Gamma)$ is now the (one-dimensional) complex space of $\Gamma$-automorphic holomorphic forms of degree $- 4$ on $\D^*$ with $\B$-norm.

All extremal Beltrami coefficients $\mu(z)$ of the maps $X_0$ onto the points $X \in \T(0, 4)$ are of the form
$\mu_t(z) = t |\psi_0(z)|/\psi_0(z)$  with $|t| < 1$ and
$$
\psi_0(z) = 1/z(z - 1)(z - \a).
$$
Hence, assuming in the contrary, that the space $\T(\Gamma)$ coincides with the disk
$\{\vp \in \B(\Gamma): \ \|\vp\|_{\B(\Gamma)} < 2\}$, one obtains biholomorphic (conformal) map of the disk $\{|t| < 1\}$ onto the ball (disk) of harmonic Beltrami coefficients
$$
\nu_{f^{\mu_t}}(z) = - \fc{1}{2} (1 - |z|^2)^2 \ S_{f^{\mu_t}}(1/\ov z), \quad |z| <  1,
$$
and by Schwarz's lemma,
$$
\|\nu_{f^{\mu_t}}\|_\iy \le |t|,
$$
which contradicts the extremality of coefficients $\mu_t$. This proves Lemma 3.

Having the existence of Schwarzians with $\|\vp\|_\B > 2$, one can straightforwardly repeat for $\T(0, 4)$ all arguments from first step, completing the proof of the theorem.

\bigskip\noindent
{\bf 4. More on geometric properties of boundary of one-dimensional space}.
As was mentioned above, the are two one-dimensional Teichm\"{u}ller spaces $\T(1, 1)$ of punctured tori and $\T(0, 4)$ of spheres with four punctures. These spaces are naturally isomorphic. Both spaces  conformally equivalent to the unit disk.

All this concerns only the hyperbolic Riemann surfaces. The assumption $2g - 2 + n > 0$ drops out the tori (equivalently, the complex elliptic curves in $\hC^2$  defined by the equation
$$
w^2 = z (z - 1) (z - \a)
$$
with $\a \in \C \setminus \{0, 1\}$). Their Teichm\"{u}ller space $\T(1, 0)$ also is one-dimensional and admits some special features; it does not be subjected to Bers' embedding.

It was established by Minsky \cite{Mi} that the boundary of either one-dimensional Teichm\"{u}ller  space $\T(1,1)$ or $\T(0, 4)$ is a Jordan curve.
A numerical investigation of points of this boundary was delivered by Porter in \cite{Po}.

Theorem 1 yields that this curve is not starlike.

\bigskip\noindent
{\bf 5. Additional remark}. Werner \cite{We} establshed that for any rectilinear or circular polygon $P$
whose sides touch a common circle and the least interior angle equals $\a \pi$, the Grunsky and Teichm\"{u}ller norms of its outer conformal mapping function $f$ have the value
$$
\vk(f) = k(f) = 1 - \a
$$
(this quantity also gives the value of the reflection coefficient across the curve $\partial $ and
the reciprocal to the first nontrivial Fredholm eigenvalue of this curve).

Lemma 2 yields that this equality is valid for all polygons (and more general quasidisks), whose mapping function $f$ satisfies $\|S_f\|_\B < 2$.

\medskip
{\small\em{ \leftline{Department of Mathematics, Bar-Ilan University, Ramat-Gan, Israel}
\leftline{and Department of Mathematics, University of Virginia, Charlottesville, VA 22904-4137, USA}}}

\end{document}